\documentclass[a4paper,10pt]{article}

\usepackage[top=4cm, bottom=3cm, left=3.2cm, right=3.2cm]{geometry}

\usepackage{Macros}

\begin{document}

\title{Modeling random traffic accidents by conservation laws}

\author{Simone G\"ottlich\footnotemark[1], \; Stephan Knapp\footnotemark[1]}
\renewcommand*{\thefootnote}{\fnsymbol{footnote}}

\footnotetext[1]{University of Mannheim, Department of Mathematics, 68131 Mannheim, Germany (goettlich@uni-mannheim.de, stknapp@mail.uni-mannheim.de).}

\date{\today}

\maketitle

\begin{abstract}
We introduce a stochastic traffic flow model to describe random traffic accidents on a single road.
The model is a piecewise deterministic process incorporating traffic accidents and is
based on a scalar conservation law with space-dependent flux function. 
Using a Lax-Friedrichs discretization, we show that the total variation is bounded in finite time and 
provide a theoretical framework to embed the stochastic process.
Additionally, a solution algorithm is introduced to also investigate the model numerically.
\end{abstract}

\noindent
{\bf AMS Classification:} 35L65, 60J25, 90B20\\ 
{\bf Keywords:} conservation laws, traffic flow, random accidents, piecewise deterministic processes 

\maketitle

\section{Introduction}
Macroscopic traffic flow models based on hyperbolic conservation laws have 
been intensively investigated during the last decades, see \cite{MR3553143,MR2328174} for an overview. 
The various research directions include theoretical and numerical investigations such for instance well-posedness~\cite{MR3447130}, 
coupled models~\cite{MR1951956}, network extensions~\cite{MR2328174,MR1338371}, optimal control~\cite{MR2164806}, or more recently, 
data-driven approaches~\cite{MR3231191} while stochastic traffic models have been less considered~\cite{JABARI201315,MR2216726}.

Typically, macroscopic traffic flow equations are either characterized by first-order models for the evolution
of the traffic density or second-order models, where an additional equation for the velocity is considered. 
So far, the modeling of traffic accidents (or incidents) has been considered in a deterministic setting~\cite{MR2842413,MR3210745,MR3116165},
queueing theory approaches~\cite{MR2474146,MR3922072} or kinetic models~\cite{MR3605555}.  
There are only a few contributions, where the presence of accidents is described by a stochastic process~\cite{MR3922072}.

Therefore, the aim of this paper is to combine the stochastic modeling of accidents with the Lighthill-Whitham-Richards (LWR) model~\cite{Lighthill1955}
of first-order type and to provide a framework that allows for theoretical and numerical studies. The idea is to include random effects
directly in the flux function such that failures depend on the current traffic density.

We assume that accidents happen at random times and have an impact on the road capacity around the accident. Based on the LWR model, we incorporate these accidents by a space-dependent flux function determining the deterministic structure between the random accidents. Obviously, the profile of the traffic density has an impact on the probability of an accident. For instance, fluctuations in the density lead to different velocities of the cars and an accident is more likely as it is the case for stationary traffic situations. 
The traffic density does not only influence the probability of an accident. It also indicates where an accident could happen as for example at the end of a traffic jam. 
In order to capture these ideas, we face two building blocks, i.e.\ the deterministic dynamics between accidents and the stochastic nature, which interrupts the deterministic flow at random times. This directly leads to the well-known piecewise deterministic processes (PDPs), see \cite{Davis1984, Jacobsen2006}. In \cite{GoettlichKnapp2019}, the latter idea has been used to incorporate random machine failures of machines based on hyperbolic dynamics, where the product density influences machine failures and vice versa. Compared to \cite{GoettlichKnapp2019}, we face different challenges here: 
First, we deal with a nonlinear dynamics with a space-dependent flux function, which does not admit total variation bounds in general and we prove under which conditions we can guarantee these bounds. Second, the position of an accident depends on the current density, which makes the modeling more involved. Additionally, the classical thinning algorithm, see \cite{Lemaire2018}, to sample the times of an accident might lead to large computational costs.

There are different works about hyperbolic equation based dynamics connected to randomness as for example random velocity fields~\cite{MR3523082,Mishra2016} and propagation of uncertainty \cite{Fjordholm2017}. However, in these works, there is no influence of the conserved quantity on the stochastic nature, i.e.\ no bi-directional relation between the deterministic and stochastic ideas.

The paper is organized as follows: in Section~\ref{sec:mode}, we present the modeling of accidents within the LWR model
and show that the total variation of the new model is bounded. Furthermore, the stochastic process is
characterized such that accident probabilities can be embedded. In Section~\ref{sec:results}, a stochastic solution algorithm based on a Lax-Friedrichs discretization is introduced 
to analyze the occurrence of traffic accidents from a numerical point of view.   

\section{Modeling of accidents}
\label{sec:mode}

We introduce how accidents as capacity drops can be incorporated into the LWR model.
As we will see, this leads to a conservation law with space-dependent flux function.
The latter equation is then extended to the possibility of a single (or multiple) random accidents.

\subsection{General setting}

Let $f \colon [0,1]\to [0,\infty)$ be a function of LWR type, i.e.
$f(\rho) = \rho(1-\rho)$ with  $f(0) = f(1) = 0$, $f^{\prime \prime} \leq \underline{c}<0$ for some $\underline{c}<0$ 
and a unique $\rho^\ast \in (0,1)$ such that $f^\prime(\rho^\ast) = 0$. 
To describe the capacities of the road, we assume a function $c_{road}\colon \R \to \R_{>0}$ and use $c_{road}(x)f(\rho)$ as space-dependent flux. 
An appropriate choice for $c_{road}$ might be piecewise constant, describing the dependency of speed limits or the number of lanes. 

We interpret an accident on a road as capacity reduction within an interval $I(p,s) \subset (p-s,p+s)$ of length $s$, where $p \in \R$ denotes the position and $s \in \R$ the size of the accident. The amount of capacity reduction is denoted by $c \in [0,c_{\text{max}}]$ with $0\leq c_{\text{max}}<1$ such that the road capacity at $p$ is given by $(1-c)c_{road}(p)$.
We denote by $x \mapsto c_a(x,p,s,c)$ the capacity function of the accident. Then, it is natural to define the space-dependent flux function 
\begin{align*}
F^{p,s,c}(x,\rho) =c_a(x,p,s,c)c_{road}(x)f(\rho).
\end{align*}

Altogether, we end up with the following Cauchy problem
 \begin{align}
 \rho_t+(F^{p,s,c}(x,\rho))_x = 0, \quad \rho(x,0) = \rho_0(x), \label{eq:IVP}
 \end{align}
 which admits a unique entropy solution, see \cite{CocliteRisebro2005} if $\TV(c_a(\cdot,p,s,c)c_{road}(\cdot))< \infty$, $\TV(\rho_0)<\infty$ and if $c_a(\cdot,p,s,c)c_{road}(\cdot)$ is differentiable with except of finitely many points. 
 Additionally, we need that 
 \begin{align*}
  \TV(\Psi(\rho(\cdot,t)))<\infty \text{ for all } t \in [0,T], \quad \Psi(\rho) = \sgn(\rho-\rho^\ast)\frac{f(\rho^\ast)-f(\rho)}{f(\rho^\ast)}.
 \end{align*}
 This does not imply $\TV(\rho)<\infty$, which we will need in the modeling of stochastic accidents later. 
However, the following lemma provides conditions on the data such that the solution to the scalar conservation law \eqref{eq:IVP} remains in $\BV(\R)$.

\begin{lemma}\label{lem:TVBound}
Let $a(x):= c_a(x,p,s,c)c_{road}(x)$ satisfy $a \in C^2(\R)\cap \TV(\R)$ and let $f$ be an LWR flux. Furthermore, we assume 
\begin{align*}
a,a^\prime,f,f^\prime \in L^\infty(\R), \quad a^\prime,a^{\prime \prime} \in L^1(\R), \quad \rho_0 \in \BV(\R).
\end{align*}
Then there exists a constant $C = C(T,\|a\|_\infty,\|a^\prime\|_\infty,\|a^{\prime \prime}\|_{1},\|f\|_\infty,\|f^\prime\|_\infty,\TV(\rho_0))$
such that the solution to \eqref{eq:IVP} satisfies $\TV(\rho(t)) \leq C$ for all $t \in [0,T]$ and $\|\rho(t)\|_\infty\leq \|\rho_0\|_\infty+T \|a^\prime\|_\infty\|f\|_\infty$. Additionally, the mapping $t \mapsto \TV(\rho(t))$ is Lipschitz continuous on $[0,T]$.
\end{lemma}
\begin{proof}
We prove the lemma by using the Lax-Friedrichs scheme given by
\begin{align} \label{eqn:LxF}
\rho_i^{j+1} = \rho_i^j-&\gl\Big(\frac{1}{2\gl}(\rho_i^j-\rho_{i+1}^j)+\frac{1}{2}(a_i f(\rho_i^j)+a_{i+1}f(\rho_{i+1}^j))\nonumber\\
&-(\frac{1}{2\gl}(\rho_{i-1}^j-\rho_{i}^j)+\frac{1}{2}(a_{i-1} f(\rho_{i-1}^j)+a_{i}f(\rho_{i}^j))\Big).
\end{align}
The convergence of the Lax-Friedrichs scheme has been studied in \cite{KARLSEN2004}, whereas in \cite{Towers2000,Towers2018} the Godunov scheme has been examined. For our purpose, the Lax-Friedrichs scheme is a suitable choice avoiding the study of various cases as needed for the Godunov scheme.
We start with the $L^\infty$ estimate followed by the $\BV$-estimate and conclude that the numerical scheme converges to the unique solution of the Cauchy problem.\\[1ex]
$\bullet$ \emph{$L^\infty$ estimate.}
Using the CFL condition $$\gl \|a\|_\infty\| f^\prime\|_\infty\leq 1,$$ we deduce that
\begin{align*}
|\rho_i^{j+1}| &= \left|\frac{\rho_{i+1}^j+\rho_{i-1}^j}{2}-\frac{\gl}{2}(a_{i+1}f(\rho_{i+1}^j)-a_{i-1}f(\rho_{i-1}^j))\right|\\
	&= \frac{1}{2}|\rho_{i+1}^j+\rho_{i-1}^j-\gl(a_{i+1}f^\prime(\xi_i)(\rho_{i+1}^j-\rho_{i-1}^j)+f(\rho_{i-1}^j)(a_{i+1}-a_{i-1}))|\\
	& \leq \frac{1}{2}\left(|\rho_{i+1}^j|(1-\gl a_{i+1}f^\prime(\xi_i))+|\rho_{i-1}^j|(1+\gl a_{i+1}f^\prime(\xi_i))\right)+\frac{\gl}{2} f(\rho_{i-1})|a^\prime(\eta_i)|2 \Delta x\\
	&\leq \|\rho^j\|_\infty+\Delta t \|f\|_\infty \|a^\prime\|_\infty.
\end{align*}
The latter implies $\|\rho^j\|_\infty \leq \|\rho_0\|_\infty+T \|f\|_\infty \|a^\prime\|_\infty$.\\[1ex]
$\bullet$ \emph{$\BV$ estimates.} Using the same arguments as in the $L^\infty$ estimates, we can estimate the spatial $\BV$ bound as follows:
\begin{align*}
\TV(\rho^{j+1}) &= \frac{1}{2}\sum_{i \in \Z} |(\rho_{i-1}^j-\rho_{i-2}^j)+\gl(a_{i-1}f(\rho_{i-1}^j)-a_{i-2}f(\rho_{i-2}^j))\\
&\quad \quad +(\rho_{i+1}^j-\rho_{i}^j)-\gl(a_{i+1}f(\rho_{i+1}^j)-a_{i}f(\rho_{i}^j))|\\
&= \frac{1}{2} \sum_{i \in \Z} |(\rho_{i-1}^j-\rho_{i-2}^j)(1+\gl a_{i-1}f^\prime(\xi^j_{i-\frac{3}{2}}))+\gl f(\rho^j_{i-2}) \Delta x a^\prime(\eta_{i-\frac{3}{2}})\\
&\quad \quad+(\rho^j_{i+1}-\rho^j_{i})(1-\gl a_{i+1}f^\prime(\xi^j_{i+\frac{1}{2}}))-\gl f(\rho^j_{i}) \Delta x a^\prime(\eta_{i+\frac{1}{2}})|\\
&\leq \frac{1}{2}\sum_{i \in \Z} |\rho_{i-1}^j-\rho_{i-2}^j||1+\gl a_{i-1}f^\prime(\xi^j_{i-\frac{3}{2}})|\\
&\quad +\frac{1}{2}\sum_{i \in \Z} |\rho_{i+1}^j-\rho_{i}^j||1-\gl a_{i+1}f^\prime(\xi^j_{i+\frac{1}{2}})|\\
&\quad +\frac{\gl \Delta x}{2}\sum_{i \in \Z} |a^\prime(\eta_{i+\frac{1}{2}})f(\rho_i^j)-a^\prime(\eta_{i-\frac{3}{2}})f(\rho_{i-2} ^j)|.
\end{align*}
 Using the CFL condition and
 \begin{align*}
  &\;|a^\prime(\eta_{i+\frac{1}{2}})f(\rho_i^j)-a^\prime(\eta_{i-\frac{3}{2}})f(\rho_{i-2} ^j)| \\
   \leq &\;\|a^\prime\|_\infty\|f^\prime\|_\infty|\rho_{i}^j-\rho_{i-2}^j|+\|f\|_\infty|a^{\prime\prime}(\tilde{\eta_i})\|3\Delta x,
 \end{align*}
yields
\begin{align*}
\TV(\rho^{j+1}) \leq  (1+\Delta t \|a^\prime\|_\infty\|f^\prime\|_\infty)\TV(\rho^j)+  \Delta t \frac{3}{2}\|f\|_\infty \|a^{\prime \prime}\|_1.
\end{align*} 
 Hence, we have
 \begin{align*}
 \TV(\rho^{j+1})& \leq e^{\|a^\prime\|_\infty\|f^\prime\|_\infty T}\TV(\rho_0)+\frac{\frac{3}{2}\|f\|_\infty \|a^{\prime \prime}\|_1}{\|a^\prime\|_\infty\|f^\prime\|_\infty }(e^{\|a^\prime\|_\infty\|f^\prime\|_\infty T}-1)\\
 &=:C_1.
 \end{align*}

Furthermore, we deduce the following bound on the time difference of the total variation
\begin{align*}
\TV(\rho^{j+m})-\TV(\rho^{j}) &= \sum_{k=0}^{m-1} (\TV(\rho^{j+k+1})-\TV(\rho^{j+k}))\\
&\leq \Delta t \sum_{k=0}^{m-1} (\|a^\prime\|_\infty\|f^\prime\|_\infty\TV(\rho^{j+k})+\frac{3}{2}\|f\|_\infty \|a^{\prime \prime}\|_1)\\
&\leq m \Delta t(C_1\|a^\prime\|_\infty\|f^\prime\|_\infty+\frac{3}{2}\|f\|_\infty \|a^{\prime \prime}\|_1).
\end{align*}
If $t =j \Delta t$ and $\tilde{t} = (j+m)\Delta t$, then $$|\TV(\rho^{j+m})-\TV(\rho^{j})| \leq \tilde{C_1} |t-\tilde{t}|.$$

 In order to use a compactness argument for the numerical scheme to converge, we need the total variation in space and time. For piecewise constant function $\rho$ it holds
 \begin{align*}
 \TV_{\R \times [0,T]}(\rho) = \sum_{j = 0}^{\frac{T}{\Delta t}} \Delta t \TV(\rho^j) +\sum_{i \in \Z} \Delta x \sum_{j = 0}^{\frac{T}{\Delta t}-1} |\rho_i^{j+1}-\rho_i^j|.
 \end{align*}
 We can directly estimate the first expression by 
 $$ \sum_{j = 0}^{\frac{T}{\Delta t}} \Delta t \TV(\rho^j) \leq T C_1.$$
 To analyze the second expression we start with
  \begin{align*}
 |\rho_i^{j+1}-\rho_i^j| &= \frac{1}{2}|(\rho_{i+1}^j-\rho_{i}^j)-\gl(a_{i+1}f(\rho_{i+1}^j)-a_i f(\rho_i^j))-(\rho_i^j-\rho_{i-1} ^j)-\gl(a_i f(\rho_i^j)-a_{i-1}f(\rho_{i-1}^j))|\\
 &=\frac{1}{2}|(\rho_{i+1}^j-\rho_{i}^j)(1-\gl a_{i+1}f^\prime(\xi_{i+\frac{1}{2}}))-(\rho_i^j-\rho_{i-1} ^j)(1+\gl a_i f^\prime(\rho_{i-\frac{1}{2}}))\\
&\quad \quad -\gl\Delta x(f(\rho_i^j)a^\prime(\eta_{i+\frac{1}{2}})+f(\rho_{i-1}^j)a^\prime(\eta_{i-\frac{1}{2}}))|\\
&\leq \frac{1}{2}|(\rho_{i+1}^j-\rho_{i}^j)|(1-\gl a_{i+1}f^\prime(\xi_{i+\frac{1}{2}}))+\frac{1}{2}|\rho_i^j-\rho_{i-1} ^j|(1+\gl a_i f^\prime(\xi_{i-\frac{1}{2}}))\\
&\quad \quad +\frac{1}{2}\gl \Delta x (f(\rho_i^j)|a^\prime(\eta_{i+\frac{1}{2}})|+f(\rho_{i-1}^j)|a^\prime(\eta_{i-\frac{1}{2}}))|),
 \end{align*}
 where we use the CFL condition and $f \geq 0$.
This leads to
\begin{align*}
\sum_{i \in \Z}|\rho_i^{j+1}-\rho_i^j| &\leq \TV(\rho^j)+\gl \Delta x\sum_{i \in \Z}f(\rho_i^j)|a^\prime(\eta_{i+\frac{1}{2}})|\\
& \leq C_1 +\gl \|f^\prime\|_\infty \|a^\prime\|_1
\end{align*}
and therefore
\begin{align*}
\sum_{i \in \Z} \Delta x \sum_{j = 0}^{\frac{T}{\Delta t}-1} |\rho_i^{j+1}-\rho_i^j| &\leq \frac{T}{\gl} C_1+ T \|f^\prime\|_\infty \|a^\prime\|_1\\
&=: C_2.
\end{align*}
Let $(\Delta t_n)_{ n \in \N}$ be a sequence, which converges to zero and $\Delta x_n = \frac{\Delta t_n}{\gl}$ be the corresponding spatial discretization, satisfying the CFL condition. The constructed sequence of piecewise constant functions $(\bar{\rho}_n)_{n \in \N}$ has a subsequence $(\bar{\rho}_{n_l})_{l \in \N}$, which converges to some $\rho \in \BV(\R \times [0,T])$ in $L^1_{loc}(\R)$ by Helly's theorem.
A Kruzkov type inequality, see \cite{KARLSEN2004}, and a Lax-Wendroff type argument show that $(\bar{\rho}_n)_{n \in \N}$ converges to a weak entropy solution, which is unique by \cite{CocliteRisebro2005}. Consequently, the limiting solution is the solution to the IVP satisfying claimed properties of the lemma.
\end{proof}
Hence, we are now able to mathematically introduce traffic accidents as partial road capacity drops via the function $a$.
 
\subsection{Random traffic accidents}

The parameters to incorporate a traffic accident in equation~\eqref{eq:IVP} are the position $p$, the size $s$ and the capacity drop $c$. From the modeling perspective the position is the first parameter to consider since there exists a dependency on the current traffic situation: if there are no cars, or cars are fully stopped by a traffic jam, we expect no accident, whereas if cars drive with high speed and the density is high at the same time, we expect a higher probability of an accident. Also, we observe accidents at the end of traffic jams. 
 To summarize, the following modeling ideas should be included:
 \begin{enumerate}
 \item a higher distance between cars at lower speed implies a lower accident probability and vice versa,
 \item a higher accident probability at increasing density (as for example tailbacks).
 \end{enumerate}
 
 \emph{Regarding 1.} The flow exactly describes the combination of density, i.e.\ car distances, and velocities such that at places where $\rho = \rho^\ast$ the probability of an accident can be assumed to be the most highest. This idea corresponds to a probability capturing random accidents caused by human failures solely (i.e. excluding tailbacks). If $v$ is uniformly bounded on $[0,1]$, the normalizing constant
 \begin{align*}
C_F := \int_\R F^{p,s,c}(x,\rho(x))dx \leq \|a\|_\infty \|v\|_\infty \int_\R \rho(x)dx 
 \end{align*}
 is finite and we can define the family of probability measures
 \begin{align}
 \mu^F_{p,s,c,\rho}(B) = \int_B \frac{1}{C_F} F^{p,s,c}(x,\rho(x))dx\label{eq:muF}
 \end{align}
 for $\rho \in \BV(\R) = \{\rho \in L^1(\R) \colon \TV(\rho)< \infty, \rho \in [0,1]\}$ and $B \in \cB(\R)$, where the latter denotes the Borel $\sigma$-algebra on $\R$. Here, we assume $\|\rho_0\|_1>0$ then it follows $C_F \neq 0$ by assumptions on $F^{p,s,c}$.
The probability measure $\mu^F_{p,s,c,\rho}$ exactly describes the probability distribution of the position of an accident caused by the flows.
 
 \emph{Regarding 2.:} In 1.\ only the information of the flow is used to specify the probability of the position of an accident. Here, we incorporate the fact that at ends of tailbacks the probability of an accident is much higher, i.e.\ if the derivative of $\rho$ is positive. Generally, for $\rho \in L^1(\R)$ we can not assign a proper derivative $D\rho$ but if $\rho \in \BV(\R)$ we can argue as follows: on the one hand, a classical derivative of $\rho \in \BV(\R)$ does not exist but on the other hand, the derivative of $\rho$ corresponds to a signed Radon measure $D\rho$ by a consequence of Riesz representation theorem. Furthermore, it holds for $ \rho \in L^1(\R)$ that
 \begin{align*}
 \TV(\rho) = \sup\left\{\int_\R \rho(x) \phi^\prime(x) dx \colon \phi \in C_c^1(\R), |\phi|\leq 1\right\} = |D\rho|,
 \end{align*}
 where $|D\rho|$ is the total variation of the measure $D\rho$ and is given by 
 \begin{align*}
 |D\rho| = D\rho^+(\R)+D\rho^-(\R).
 \end{align*}
 In the latter equation we used the Hahn decomposition, i.e.\ there exists a measurable set $\tilde{B}\in \cB(\R)$ such that $\rho^+(B) = D\rho(B\cap E)\geq 0$ and $\rho^-(B) = -D\rho(B\cap (\R\setminus E))\geq 0$ satisfy $D\rho(B) = D\rho^+(B)-D\rho^-(B)$ for every $B \in \cB(\R)$. For further details, we refer the reader to \cite{Bauer2001,EvansGariepy2015,Halmos1978,Walter1987}.
 
 A natural probability measure for $\rho \in \BV(\R)$ to describe positions of potential accidents caused by increasing densities is then given by
 \begin{align*}
 \mu^D_\rho(B) = \frac{D\rho^+(B)}{D\rho^+(\R)},
 \end{align*}
 for every $B \in \cB(\R)$ provided $D\rho^+(\R)>0$.
 
 Summarizing, we define 
 \begin{align} \label{eqn:beta}
 \mu^{pos}_{p,s,c,\rho}(B) = 
 \beta \mu^F_{p,s,c,\rho}(B)+(1-\beta) \mu^D_\rho(B) 
 \end{align}
for some fixed $\beta \in [0,1]$. That means, if $\beta = 1$, the influence of increasing densities is neglected (end of tailbacks) and if $\beta = 0$, only the latter effect is incorporated. The case $D\rho^+(\R) = 0$ means that there is no increasing part in the function $\rho$, which implies together with $\rho \in L^1(\R)$ and $\TV(\rho)<\infty$ that only $\rho = 0$ can fulfill $D\rho^+(\R) = 0$.

We only have discussed the probability distribution for the position $p$ of the accidents so far. We assume that the size $s$ follows the probability distribution $\mu^{size}$ on $(\R,\cB(\R))$ and the capacity reduction $c$ follows $\mu^{cap}$ on $([0,1),\cB([0,1)))$. In a natural way, we collect the details using the product space
\begin{align*}
E = \R \times \R \times [0,1) \times \BV(\R)
\end{align*}
with norm $$\|y\|_E = |p|+|s|+|c|+\|\rho\|_{L^1(\R)}+\TV(\rho),$$
for $y = (p,s,c,\rho) \in E$ to define a Banach space $E$. Furthermore, we denote by $\cE = \sigma(E)$ the smallest $\sigma$-algebra generated by the open sets induced by the norm $\|\cdot\|_E$. 
Finally, we define for every $y \in E$ and every $B \in \cE$ the product measure
\begin{align*}
\eta(y,B) = \mu^{pos}_{y} \otimes \mu^{size} \otimes \mu^{cap} \otimes \epsilon_\rho (B),
\end{align*}
where $\epsilon_z$ is the Dirac measure with unit mass in $z$. 
Since $\eta(y,B)$ describes the transition from no accident to one accident, we expect $\eta$ to be a kernel as the following lemma shows. 
\begin{lemma}\label{lem:measurable}
Let $(p,s,c) \mapsto \int_\R c_a(x,p,s,c)dx$ be continuous. Then $\eta$ defines a Markovian kernel on $(E,\cE)$, which additionally satisfies $\eta(y,\{y\}) = 0$ for every $y \in E$ if either $\mu^{size}(\{s\}) = 0$ for all $s$ or $\mu^{cap}(\{c\}) = 0$ for all $c \in [0,1)$.
\end{lemma}
\begin{proof}
Let $y \in E$, then $B \mapsto \eta(y,B) \geq 0$ is a measure and also $\eta(y,E) = 1$ by construction. Given a set $B \in \cE$, the mapping $y \mapsto \eta(y,B)$ is measurable if $y \mapsto \mu_y^{pos}$ is measurable since $\epsilon_\rho$ is measurable in $\rho$. 
We have $\eta(y,\{y\}) = \mu_\rho^p(\{p\}) \mu^s(\{s\}) \mu^c (\{c\}) \gep_\rho(\{\rho\}) = 0$.

It remains to show that $\rho \mapsto \mu^p_\rho$ is measurable. For every $B \in \cB(\R)$ one verifies for $\rho \neq 0$ that
\begin{align*}
0 < \mu^F_{p,s,c,\rho}(B) \leq 1, \quad 0 < \mu^D_{\rho}(B) \leq 1.
\end{align*}
Take $y = (p,s,c,\rho)$, $\tilde{y} = (\tilde{p},\tilde{s},\tilde{c},\tilde{\rho}) \in E$, satisfying $\rho$, $\tilde{\rho} \neq 0$. We deduce
\begin{align*}
|\mu_y^F(B)-\mu_{\tilde{y}}^F(B)| &\leq \frac{2}{\|F^{p,s,c}(\cdot,\rho(\cdot))\|_1}\left(\|F^{p,s,c}(\cdot,\rho(\cdot))-F^{\tilde{p},\tilde{s},\tilde{c}}(\cdot,\tilde{\rho}(\cdot))\|_1\right)\\
& \leq \frac{2}{\|F^{p,s,c}(\cdot,\rho(\cdot))\|_1}(\|c_{road}\|_\infty\|v\|_{\infty} \|\rho\|_1 \|c_a(\cdot,p,s,c)-c_a(\cdot,\tilde{p},\tilde{s},\tilde{c})\|_1\\
&\quad \quad +\|f^\prime\|_\infty \|c_{road}\|_\infty \|\rho-\tilde{\rho}\|_1).
\end{align*}
We also have
\begin{align*}
|\mu_\rho^D(B)-\mu_{\tilde{\rho}}^D(B)|& \leq \frac{1}{D\rho^+(\R)}(|D\rho^+(B)-D\tilde{\rho}^+(B)|+|D\rho^+(\R)-D\tilde{\rho}^+(\R)|)\\
&\leq  \frac{1}{D\rho^+(\R)}\TV(\rho-\tilde{\rho}).
\end{align*}
Hence, the mapping $y \mapsto \mu^{pos}_y(B)$ is continuous and therefore measurable.
\end{proof}

So far, we only have specified the probability distribution of a jump in the case that a jump occurs. To construct the time of a jump, or accident, we additionally need information about how likely a jump at time $t$ is. This can be done with rate functions and is based on the ideas of a marked point process, or, deterministic Markov processes, see \cite{Davis1984,Jacobsen2006}.

A possible choice for a rate function $\psi \colon E  \to (0,\infty)$ is given by 
\begin{align*}
\psi(y) = \gl^FC_F(\rho)+\gl^D D\rho^+(\R),
\end{align*}
where $\gl^F,\gl^D>0$ scale the influence of accidents caused by high fluxes and ends of tailbacks, respectively. For fixed $y = (p,s,c,\rho) \in E$, the rate $\psi(y)$ is finite. More precisely, if $\bar{\rho}(x,t)$ is a weak entropy solution to the IVP \eqref{eq:IVP}, then for $a(x) = c_a(x,p,s,c)c_{road}(x)$ it holds that
\begin{align*}
& \;\gl^FC_F(\bar{\rho}(t))+\gl^D D(\bar{\rho}(t))^+(\R)  \\
 \leq &\; \gl^F \|a\|_\infty \|v\|_\infty \int_\R \rho_0(x)dx +\gl^D \TV(\rho(t)) \\
 \leq &\; \gl^F \|a\|_\infty \|v\|_\infty \|\rho\|_1 +\gl^D C(T,\|a\|_\infty,\|a^\prime\|_\infty,\|a^{\prime \prime}\|_{1},\|f\|_\infty,\|f^\prime\|_\infty,\TV(\rho))\\
 =:&\; \bar{\gl}(y).
 \end{align*}
 We have to keep in mind that for $y = (p,s,c,\rho) \in E$ the values $\|a\|_\infty,\|a^\prime\|_\infty,\|a^{\prime \prime}\|_{1}$ might differ. 
We know that  $\|a\|_\infty = 1$ and $a^{\prime \prime} = 0$ for all $x \in \R \setminus I(p,s)$ by assumption. Hence, $\|a^{\prime \prime}\|_1 \leq |I(p,s)|\||a^{\prime \prime}\|_\infty$.
Therefore, we assume $a \in C^2(\R)$, cf. Lemma~\ref{lem:TVBound}.
 
Let $\phi \colon E \to E$ be the deterministic evolution, i.e.\ 
\begin{align*}
\phi_t((p_0,s_0,c_0,\rho_0)) = (p_0,s_0,c_0,\rho(t)),
\end{align*}
where $\rho(t)$ is the unique weak entropy solution to the IVP \eqref{eq:IVP} with initial datum $\rho_0$ and the parameters $p_0,s_0,c^{\max}  = c_0$.

Let $(U_i, i \in \N)$ be a sequence of independent and identically distributed (i.i.d) random variables on some probability space \OAP\ each having a uniform distribution on $[0,1]$. Furthermore, let $(\xi_i, i \in \N)$ be a sequence of i.i.d exponentially distributed random variables on the same probability space \OAP\ and independent of $(U_i, i \in \N)$ and choose $t_n \in [0,T]$, $y_n \in E$. The following thinning algorithm produces the next jump time $T_{n+1}$ and corresponding post jump location $Y_{n+1}$.

\begin{algorithm}[H]
\caption{Thinning algorithm}
\label{alg:thinning_one_queue_proc}
\begin{algorithmic}
\STATE $i = 1$
\STATE $s_i = t_n+\xi_i$
\WHILE {$U_i > \psi(\phi_{t_n s_i}(y_n)) \cdot (\overline{\lambda})^{-1}$ and $s_i<T$}
	\STATE $s_{i+1} = s_i + \xi_i$
	\STATE $i = i+1$
\ENDWHILE
\STATE $T_{n+1} = s_i$
\STATE Generate $Y_{n+1} \sim \eta(\phi_{t_n s_i}(y_n),\cdot)$
\end{algorithmic}
\end{algorithm}
One can show, see \cite{GoettlichKnapp2019}, that
\begin{align}
P(T_{n+1}\leq t) &= 1-e^{-\int_{t_n}^t \psi(\phi_{\tau-t_n}(y_n))d\tau},\notag \\
P(Y_{n+1} \in B|T_{n+1} = t) &= \eta(\phi_{t-t_n}(y_n),B) \label{eq:ThinningJumpDist}
\end{align}
for $t \geq t_n$ and $B \in \cE$.

We set  $T_0 = 0$ and $Y_0 = (p_0,s_0,c_0,\rho_0) \in E$ and apply the thinning algorithm iteratively. In every iteration we obtain a new upper bound $\bar{\gl}$ on the rates, which might increase but stays finite for finitely many iterations. Let denote $((T_n,Y_n),n\in \N_0)$ the constructed jump times and post-jump locations, then we define the piecewise deterministic process (PDP) $(X(t), t \in [0,T])$ as
\begin{align*}
X(t) = Y_n \Leftrightarrow t \in [T_n,T_{n+1}).
\end{align*}

\begin{remark}\hspace{0mm}
\begin{enumerate}
\item The total variation bound on the solution is quite pessimistic for reasonable initial datum. 
\item The total variation bound can be very large in small time intervals and the Algorithm \ref{alg:thinning_one_queue_proc} can not be used efficiently to simulate the model.
\item We expect $X$ being a Markov process but standard results, see \cite{Jacobsen2006} can not be applied since $\BV$ is no Borel space and the existence of regular conditional distributions is not guaranteed. 
\end{enumerate}
\end{remark}

\emph{Multiple accidents on roads.}
In order to implement multiple accidents in the model, we label accidents and extend the state space as follows:
\begin{itemize}
\item positions are now given by $\vec{p} \in \R^\N$,
\item sizes of the accidents are $\vec{s} \in \R^\N$,
\item capacity reductions $\vec{c} \in [0,1)^\N$
\end{itemize}
and set
\begin{align*}
E = \R^\N \times \R^\N \times [0,1)^\N \times \BV(\R)
\end{align*}
with the norm
\begin{align*}
\|y\|_E = \|\vec{p}\|_{l^1}+\|\vec{s}\|_{l^1}+\|\vec{c}\|_{l^1}+\|\rho\|_{L^1(\R)}+\TV(\rho).
\end{align*}
Let $\gl_A>0$ be the rate of an accident and $\gl_R>0$ be the rate of resolving an accident. We define $m(\vec{c}) = \min \{ i \colon c_i = 0\}$ 
and $\pi_i(z,\vec{v}) = (v_1,\dots,v_{i-1},z,v_{i+1},\dots) \in \R^\N$. A natural choice for the jump distribution is then given by
\begin{align}
\eta(y,B) &= \frac{1}{\gl_R \sum_{i \in \N}\Ind_{c_i>0}+\gl_A}\Big[  \gl_R \sum_{i \in \N} \Ind_{c_i>0} \gep_{(\vec{p},\vec{s},\pi_i(0,\vec{c}),\rho)}(B) \notag\\
  &\quad + \gl_A \int_{\R^2\times [0,1)}  \gep_{(\pi_{m(\vec{c})}(\tilde{p},\vec{p}),\pi_{m(\vec{c})}(\tilde{s},\vec{s}),\pi_{m(\vec{c})}(\tilde{c},\vec{c}),\rho)}(B) \mu^{pos}_y \otimes \mu^{size} \otimes \mu^{cap}(d(\tilde{p},\tilde{s},\tilde{c})) \Big]. \label{eq:EtaMultiAccidents}
\end{align}
Here, $\mu^{pos}_y = \beta \mu^F_{y}+(1-\beta) \mu^D_\rho $, where $\mu^F_{y}(B) = \int_B \frac{1}{C_F} F^{\vec{p},\vec{s},\vec{c}}(x,\rho(x))dx$ and $F^{\vec{p},\vec{s},\vec{c}}(x,\rho) = c_{road}(x)f(\rho)\prod_{i \in \N}c_a(x,p_i,s_i,c_i)$.
The sum $N(\vec{c}) =  \sum_{i \in \N}\Ind_{c_i>0}$ corresponds to the number of accidents and we see that $B \mapsto \eta(y,B)$ is a probability measure. Since $\pi_i$ and $m$ are measurable functions, the mapping $y \mapsto \eta(y,B)$ is measurable if again $y \mapsto \mu^{pos}_y$ is measurable, see Lemma \ref{lem:measurable}. Since $\gl_A$ corresponds to the rate of an accident, we choose again $$\gl_A(y) =  \gl^FC_F(\rho)+\gl^D D\rho^+(\R)$$
and
\begin{align*}
\psi(y) = \gl^FC_F(\rho)+\gl^D D\rho^+(\R) + \gl_R \sum_{i \in \N}\Ind_{c_i>0}.
\end{align*}
The upper bound on the rate function is now given by 
\begin{align*}
\psi(y) \leq \gl^F \|a\|_\infty \|v\|_\infty \|\rho_0\|_1+\gl^D C(T,\|a\|_\infty,\|a^\prime\|_\infty,\|a^{\prime \prime}\|_{1},\|f\|_\infty,\|f^\prime\|_\infty,\TV(\rho_0))+\gl_R N(\vec{c}),
\end{align*}
where $a(x) = c_{road}(x)\prod_{i \in \N}c_a(x,p_i,s_i,c_i)$, $y = (\vec{p},\vec{s},\vec{c},\rho(t))$ and $\rho(t)$ is the unique weak entropy solution to \eqref{eq:IVP}.

We explain the choice of \eqref{eq:EtaMultiAccidents} by the following example.
We consider two accidents  with capacity reduction $\vec{c} = (0.5,0,0.5,0\dots)$, i.e.\ $N(\vec{c}) = 2$ and $m(\vec{c}) = 2$. We set $B_{\vec{p}} = \R \times B_p \times \R \times \cdots \in \sigma(\R^\N)$, $B_{\vec{s}} = \R \times B_s \times \R \times \cdots \in \sigma(\R^\N)$ and $B_{\vec{c}} = B_{c_1}\times B_{c_2} \times B_{c_3} \times \R \times \cdots \in \sigma([0,1)^\N)$. Then, we set $B = B_{\vec{p}}\times B_{\vec{s}} \times B_{\vec{c}} \times \BV(\R)$ and obtain
\begin{align*}
\eta(y,B) = \frac{1}{2 \gl_R +\gl_A}[\gl_R (\gep_0(B_{c_1})+\gep_0(B_{c_3}))+\gl_A\mu^{p}_\rho \otimes \mu^s \otimes \mu^c(B_p \times B_s \times B_{c_2}) ].
\end{align*}
This implies that the probability of resolving the first accident and no new accident, i.e.\ $B_{c_1} = \{0\}$, $B_{c_3} = B_{c_2} = \emptyset$, is given by
\begin{align*}
\eta(y,B) = \frac{\gl_R}{2\gl_R+\gl_A}.
\end{align*}
In the same manner we obtain the probability of having a new accident somewhere with some size and no repairs, i.e.\ $B_{c_1} = B_{c_3} = \emptyset$, $B_{c_2} = B_p = \R$ and $B_s = [0,1)$, 
\begin{align*}
\eta(y,B) = \frac{\gl_A}{2\gl_R +\gl_A}.
\end{align*}
Hence, if $\gl_A = \gl_R$, the probabilities are equal with value $\frac{1}{3}$. 

\section{Numerical treatment and computational results}
\label{sec:results}

The Cauchy problem \eqref{eq:IVP} is numerically solved using the Lax-Friedrichs scheme 
with a temporal step size $\Delta t >0$ and a fixed relation $\frac{\Delta t}{\Delta x}$ such that the scheme converges to the weak entropy solution $\rho$ of the Cauchy problem,
cf.\ Lemma \ref{lem:TVBound}.
%
%
We denote by 
$${\rho}_i^0 = \frac{1}{\Delta x} \int_{x_{i-\nicefrac{2}{1}}}^{x_{i+\nicefrac{1}{2}}} \rho_0(x) dx$$
the cell means of the initial datum $\rho_0$ for $X_i = i \Delta x$ and $i \in \Z$. 
%
%

Since the position, size and capacity reduction stays constant between the jumps, we define the discrete deterministic dynamics as
\begin{align*}
\phi^{\Delta t}_t (p_0,s_0, c_0,{\rho}_0) = (p_0,s_0,c_0,{\rho}(t)),
\end{align*}
where ${\rho}_0$ is a piecewise constant function on $[x_{i-\nicefrac{1}{2}},x_{i+\nicefrac{1}{2}})$
given by the cell means ${\rho}_i^0$. Further, ${\rho}(t)$ is the piecewise constant function given by the numerical scheme with step size $\Delta t$ and a possibly smaller last step size to reach exactly $t$.

Then, we then approximate $\mu^F_{y}$ by 
\begin{align*}
\bar{\mu}^F_{\bar{y}} (B) = \sum_{i \in \Z} \frac{F^{\vec{p},\vec{s},\vec{c}}(x_i,{\rho}_i)}{\bar{C}_F} \int_B \Ind_{[x_{i-\nicefrac{1}{2}},x_{i+\nicefrac{1}{2}})}(x)dx
\end{align*}
and 
\begin{align*}
{\bar{C}_F}=\sum_{i \in \Z}F^{\bar{p},\bar{s},\bar{c}}(x_i,\rho_i) \Delta x.
\end{align*}

Thanks to the piecewise constant cell averages, we enjoy an explicit representation of  $D\rho^+$ as
\begin{align*}
D{\rho}^+ = \sum_{i \in \Z} ({\rho}_i-{\rho}_{i-1})_+ \epsilon_{x_{i-\nicefrac{1}{2}}}\quad \text{ and }\quad  \bar{\mu}^D_{\rho} = \frac{D{\rho}^+}{D{\rho}^+(\R)}.
\end{align*}
The discretized version of the rate function $\psi(y)$ is then given  by
\begin{align*}
\bar{\psi}(y) =  \gl^F \bar{C}_F+\gl^D D{\rho}^+(\R) + \gl_R \sum_{i \in \N}\Ind_{c_i>0}.
\end{align*}
In order to use Algorithm \ref{alg:thinning_one_queue_proc}, we need a uniform upper bound on $\bar{\psi}$ which will depend on the number of accidents and grows exponentially due to the total variation bound in Lemma \ref{lem:TVBound}. In \cite{Lemaire2018}, less restrictive bounds have been used to define an appropriate algorithm but the bounds propsed will also depend on the exponential growth of the estimation of the total variation. We will introduce an approximate scheme, where the jump times are not simulated exactly in the following. The idea is based on the simulation algorithm introduced in \cite{DegondRinghofer2007}, where an algorithm has been proposed to approximate a continuous-time Markov Chain. 

The probability that an accident occurs at a time $T_{n+1}$, which is before $T_n+\Delta t$ is given by 
\begin{align}
P(T_{n+1} \leq T_n+\Delta t) = 1-e^{\int_{T_n}^{T_n+\Delta t} \psi(\phi_{\tau-T_n}(Y_n)d\tau)} = \Delta t \psi(Y_n) + \lano(\Delta t) \label{eq:AsymptoticProb1}
\end{align}
as $\Delta t \to 0$.
This is true since $t \mapsto \psi(\phi_t(Y_n))$ is Lipschitz continuous by using Lemma \ref{lem:TVBound} and the properties of $C_F$, i.e.
\begin{align*}
|\psi(\phi_t(Y_n))-\psi(\phi_{\tilde{t}}(Y_n))| \leq C(\|\rho(t)-\rho(\tilde{t})\|_1+\TV(\rho(t))-\TV(\rho(\tilde{t}))) \leq \tilde{C} |t-\tilde{t}|.
\end{align*}

Equation \eqref{eq:AsymptoticProb1} motivates the following algorithm to approximate the next jump time $T^a_{n+1}$.
\begin{algorithm}[H]
\caption{Approximate algorithm jump times}
\label{alg:approx_jumps}
\begin{algorithmic}
\STATE $i = 1$, $y = Y_n$, $t_{loc} = T_n$, $\Delta t =  \min\{\Delta t_{ref},\frac{\varrho}{\psi(y)},T-t_{loc}\}$
\WHILE {$U_i > \Delta t \psi(y)$ \AND $t_{loc}<T$}
	\STATE $t_{loc}:= t_{loc}+ \Delta t$
	\STATE $y := \phi_{\Delta t}(y)$
	\STATE $\Delta t:= \min\{\Delta t_{ref},\frac{\varrho}{\psi(y)},T-t_{loc}\}$
	\STATE $i := i+1$
\ENDWHILE
\STATE $T^a_{n+1} = t_{loc}+\Delta t$
\STATE $y := \phi_{\Delta t}(y)$
\STATE Generate $Y_{n+1} \sim \eta(y,\cdot)$
\end{algorithmic}
\end{algorithm}
The parameters $\varrho \in (0,1]$, $\Delta t_{ref} >0$ are user-defined and $(U_i, i \in \N)$ is a sequence of i.i.d.\ uniformly distributed random variables. The parameter $\Delta t_{ref}$ allows to control the accuracy of the algorithm as the reference step size and $\varrho$ is the acceptance ratio in the case that $\Delta t_{ref}$ and $T$ are large. We see that Algorithm \ref{alg:approx_jumps} uses an adaptive step size, where the adaptivity is incorporated by the current value of the rate function $\psi(y)$. We do not need any uniform bound, which is the obvious advantage and reduces the computational costs. Note that the exact solution operator $\phi$ has to be replaced by the discrete one in numerical implementations.

It remains to introduce the simulation procedure in the case that an accident happens or an accident does not cause capacity drop anymore, i.e.\ the simulation of $\eta(y,\cdot)$. The highest index $i$, where $c_i>0$ and $c_i = 0$ corresponds exactly to $N(\vec{c})$ by construction if we start with $c_{j}>0$ for $j = 1,\dots N( \vec{c})$ and $c_j = 0$ for $j > N( \vec{c})$.
One can use the well-known composition method, i.e.\ the distribution is a weighted sum of distributions, and we obtain the following procedure:
\begin{enumerate}
\item Choose whether an accident happens $Z_1 = 1$ or an accident is resolved $Z_1 = 0$ by a Bernoulli distributed random variable with $P(Z_1 = 1) = \frac{\gl_A}{\gl_R N(\vec{c})+\gl_A}$.
\item \begin{itemize}
\item Case $Z_1 = 1$: Choose independently a position $p_{N(\vec{c})+1}$ according to the law $\mu_y^{pos}$, a size $s_{N(\vec{c})+1}$ according to $\mu^{size}$ and $c_{N(\vec{c})+1} \sim\mu^{cap}$ the corresponding capacity drop.
\item Case $Z_1 = 0$: Choose a uniformly distributed index on $\{1,\ldots,N(\vec{c})\}$ to indicate which accident got removed. 
\end{itemize} 
\end{enumerate}
Simulating the new position is straightforward since $i$ is picked according to $$\sum_{i \in \Z} \frac{F(x_i,{\rho}_i)}{\bar{C}_F} \epsilon_{x_i}$$ and then the position within cell $i$ as a uniform distribution on $[x_{i-\nicefrac{1}{2}},x_{i+\nicefrac{1}{2}})$.

\subsection{Simulation results}

We assume a bounded road $[-L,L] \subset \R$ in the following with periodic boundary conditions $\rho(-L,t) = \rho(L,t)$ for~\eqref{eq:IVP} to avoid difficulties with boundary treatment. 
We assume possibly different road capacities on $[-L,L]$, i.e.\ let $$\tilde{c}_{road}(x) = \sum_{m = 0}^{M-1} c_{m,road}\Ind_{[x_m,x_{m+1})}$$ for 
$-L = x_0<x_1< \cdots x_M = L$ with $c_m \geq \underline{c}$ for $m = 0,\dots M-1$ and $c_0 = c_{M-1}$. The latter condition avoids a discontinuity for the periodic boundary conditions and implies that cars leaving at $x=L$ enter in the same manner at $x=-L$ again.
Since we need enough regularity on $c_{road}$ to apply the total variation bound on the solution of \eqref{eq:IVP}, we use a mollifier $M_\epsilon$ with support $[-\gep,\gep]$ and $\int_\R M_\gep(x) dx = 1$. Then, $c_{road}(x) = \tilde{c}_{road} \ast M_{\gep}(x) = \int_\R \tilde{c}_{road}(y) M_{\gep}(x-y)dy \in C^\infty$ and $|\operatorname{supp}(c_{road}^{\prime \prime})| \leq 2 \gep M$. 

We use the same ideas for the capacity reduction and define $\tilde{c}_a(x,p,s,c) = 1-c\Ind_{(p-\frac{s}{2},p+\frac{s}{2})}(x)$ for $p \in [-L,L]$, $s \in (-L,L)$ and $c \in [0,1-c_{min}]$.
By defining $c_a(x,p,s,c)=\tilde{c}_a(\cdot,p,s,c)\ast M_\gep(x)$, and using $a(x) = c_{road}(x)\prod_{i \in \N}c_a(x,p_i,s_i,c_i)$, we deduce 
\begin{align*}
a,a^\prime \in L^\infty(\R), \quad a^\prime, a^{\prime \prime} \in L^1(\R)
\end{align*}
as required.

\begin{remark}
We face only finitely many accidents $P$-a.s.\ such that the infinite product in $a$ can be represented by a finite product. Therefore, the differentiation of $a$ can be understood in the classical sense.
\end{remark}

The first example is devoted to the understanding of the dynamics of the LWR model with accidents derived in the previous sections. We are interested whether the modeling ideas can be also observed in computational experiments. The data we use is as follows: a time horizon $T = 60$, a spatial discretization $\Delta x = \frac{1}{50}$ of $[-10,10]$ and $\Delta t	_{ref} = \frac{1}{20}$. The initial density is chosen constant as $\rho_0 (x) = 0.4$ and the LWR flux is given by $f(\rho) = \rho(1-\rho)$. We assume a road capacity given by the non-smooth version as 
\begin{align*}
\tilde{c}_{road}(x)  = 7-2\Ind_{[0,5]}(x),
\end{align*}
which implies a capacity reduction on $[0,5]$ caused by e.g.\ roads under constructions. To incorporate capacity drops caused by accidents, we use the function
\begin{align*}
\tilde{c}_a(x,p,s,c) = 1-c\Ind_{[p-\frac{s}{2},p+\frac{s}{2}]}(x).
\end{align*}
In numerical investigations, we have recovered that smoothing the latter functions does not significantly change the results for a fixed spatial step size $\Delta x$ and $\epsilon< \frac{\Delta x}{2}$, which reduces the computational costs significantly. For the stochastic part, we use $\gl_R = \frac{1}{2}$, $\gl_D = \frac{1}{10}$, $\gl_F = \frac{1}{105}$ and assume
\begin{align}
\mu^s = \frac{1}{0.8}\Ind _{[0.2,1]}(x)dx, \quad \mu^c = \frac{1}{2}(\varepsilon_{0.5}+\varepsilon_{0.99}), \label{eq:muEx1}
\end{align}    
as well as $\varrho = 1$.

\textit{A first insight into the behavior of the model.}
Having all the parameters at hand, except $\beta$ from equation~\eqref{eqn:beta}, we can get first insights into the behavior of the model using 
numerical simulations for varying $\beta.$
The latter parameter describes the influence of the current flux on the position of possible accidents, see \eqref{eq:muF}. 

\begin{figure}[H]
\subfigure[$t = 4.9$: first accident.]{
\includegraphics[width = 0.45\textwidth]{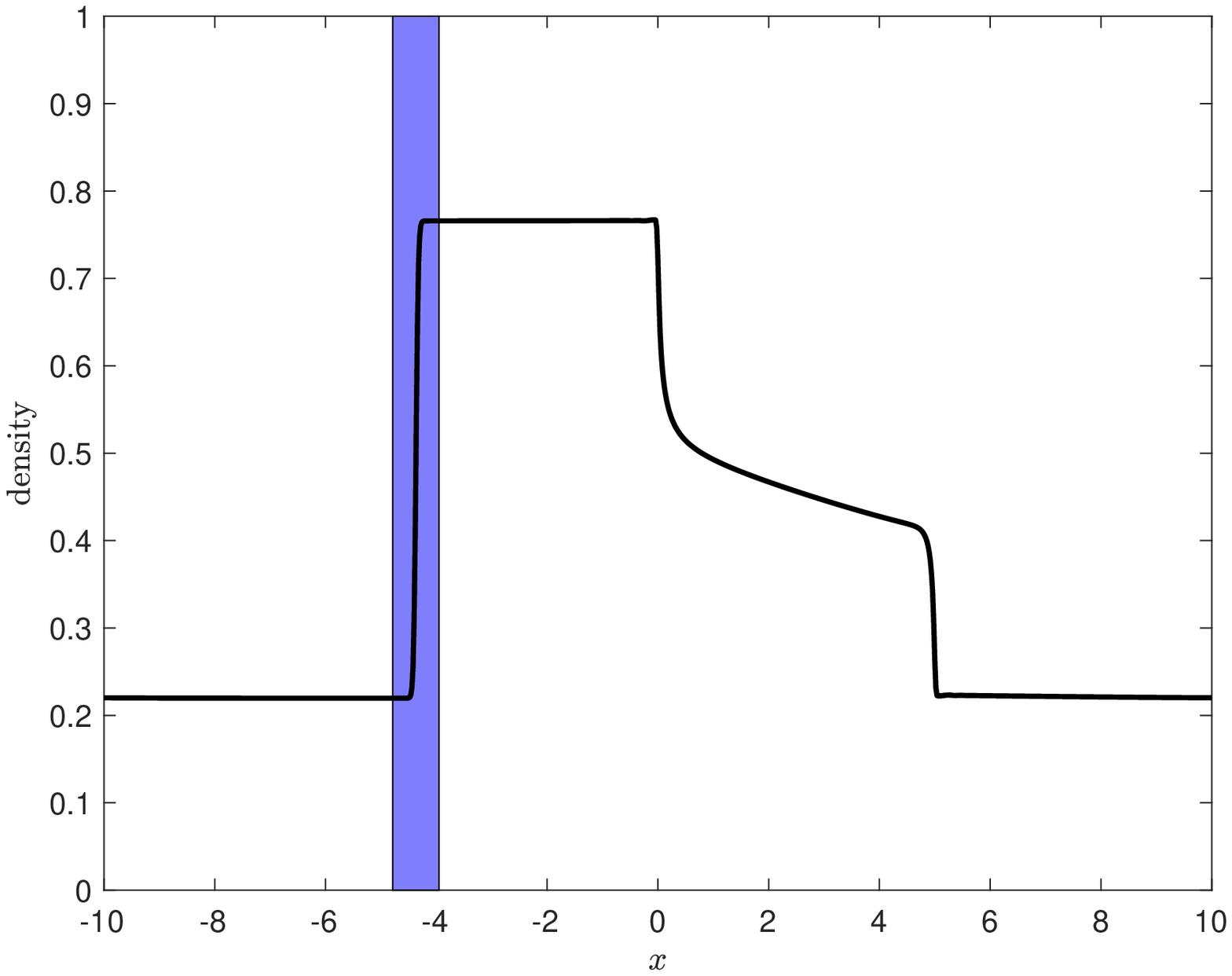}
\label{fig:Ex2SamplePathBeta0a}
}
\hfill
\subfigure[$t = 5.85$: first accident removed.]{
\includegraphics[width = 0.45\textwidth]{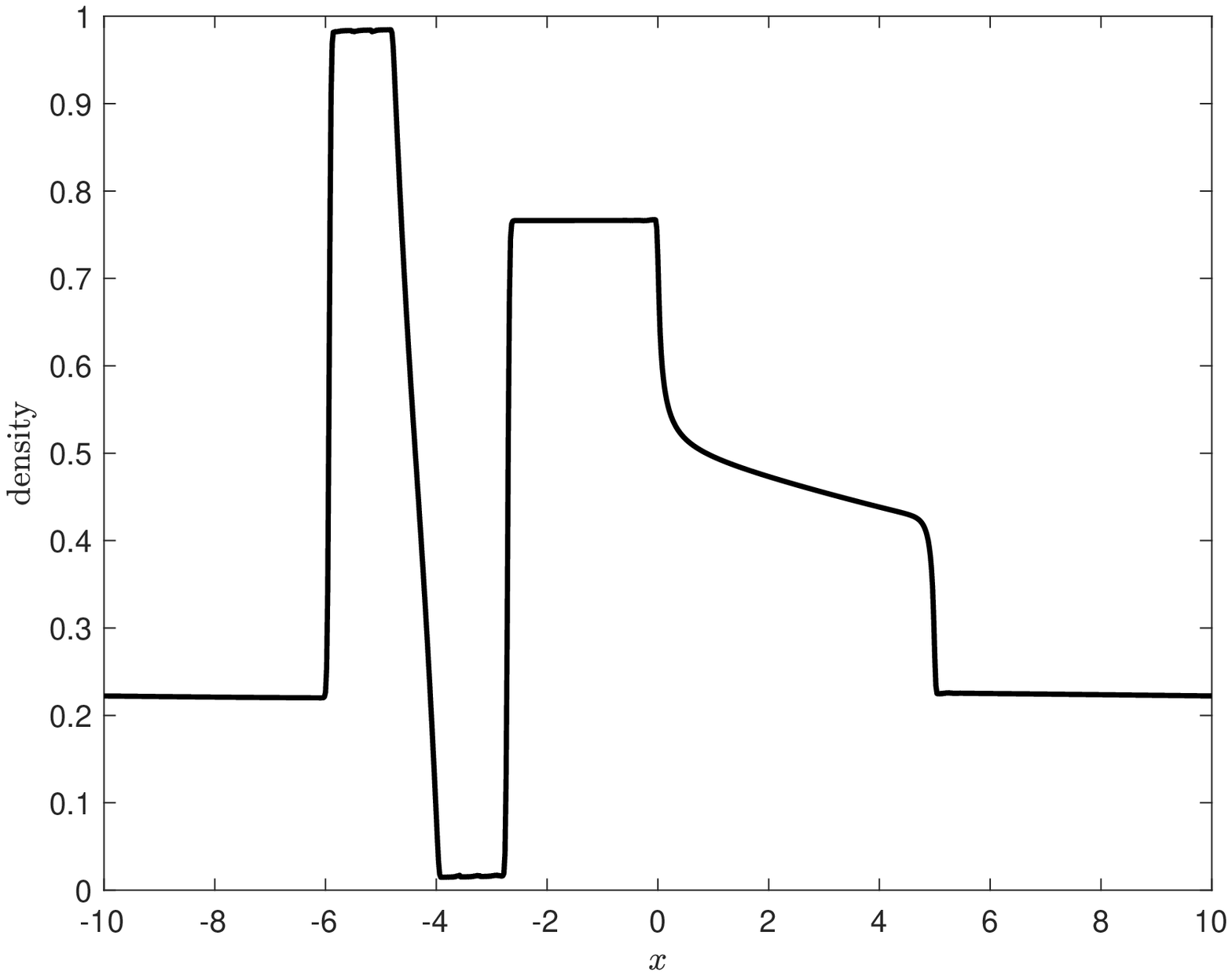}
\label{fig:Ex2SamplePathBeta0b}
}
\end{figure}

\begin{figure}[H]
\subfigure[$t = 22$: second accident.]{
\includegraphics[width = 0.45\textwidth]{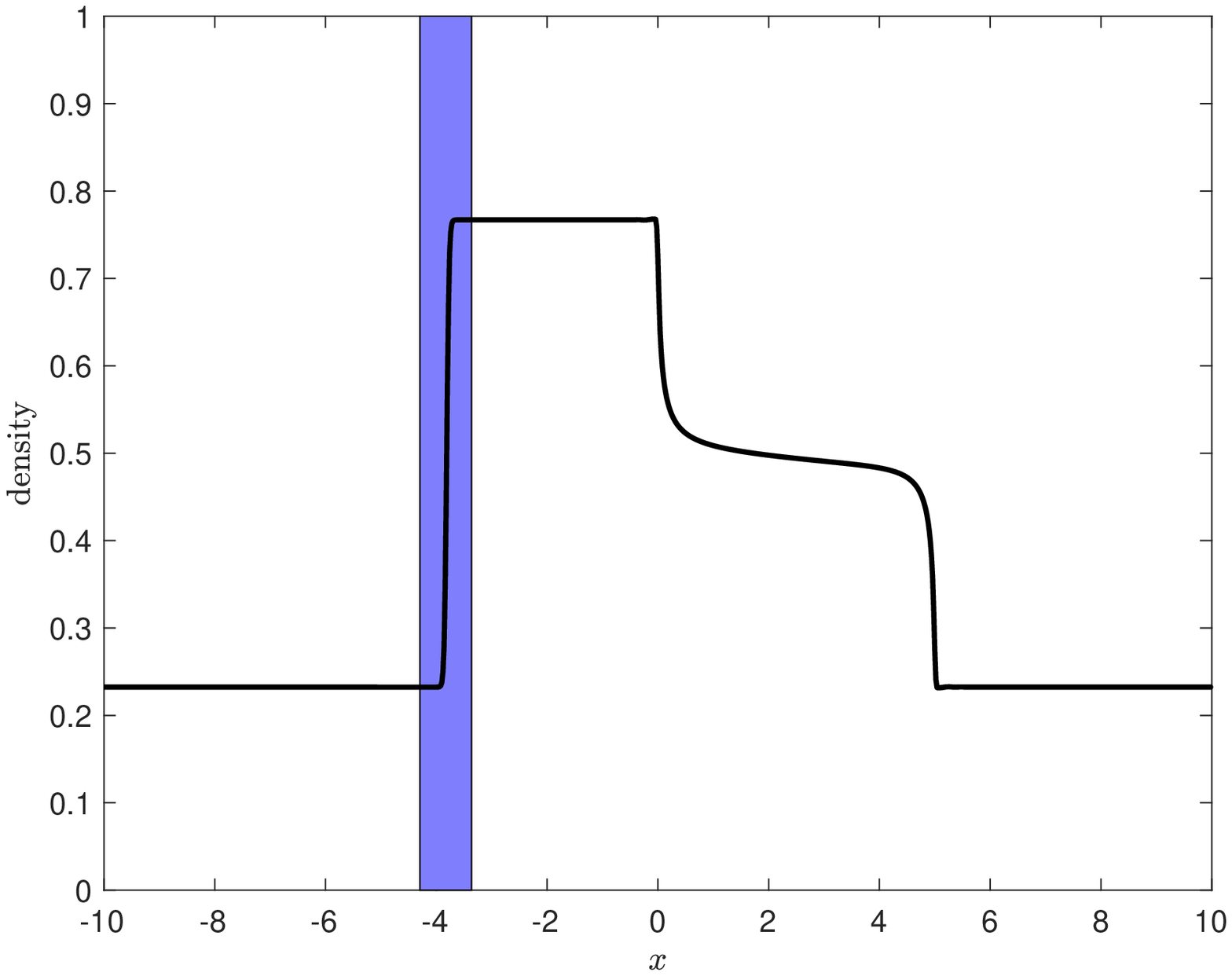}
\label{fig:Ex2SamplePathBeta0c}
}\hfill
\subfigure[$t = 26.35$: third accident.]{
\includegraphics[width = 0.45\textwidth]{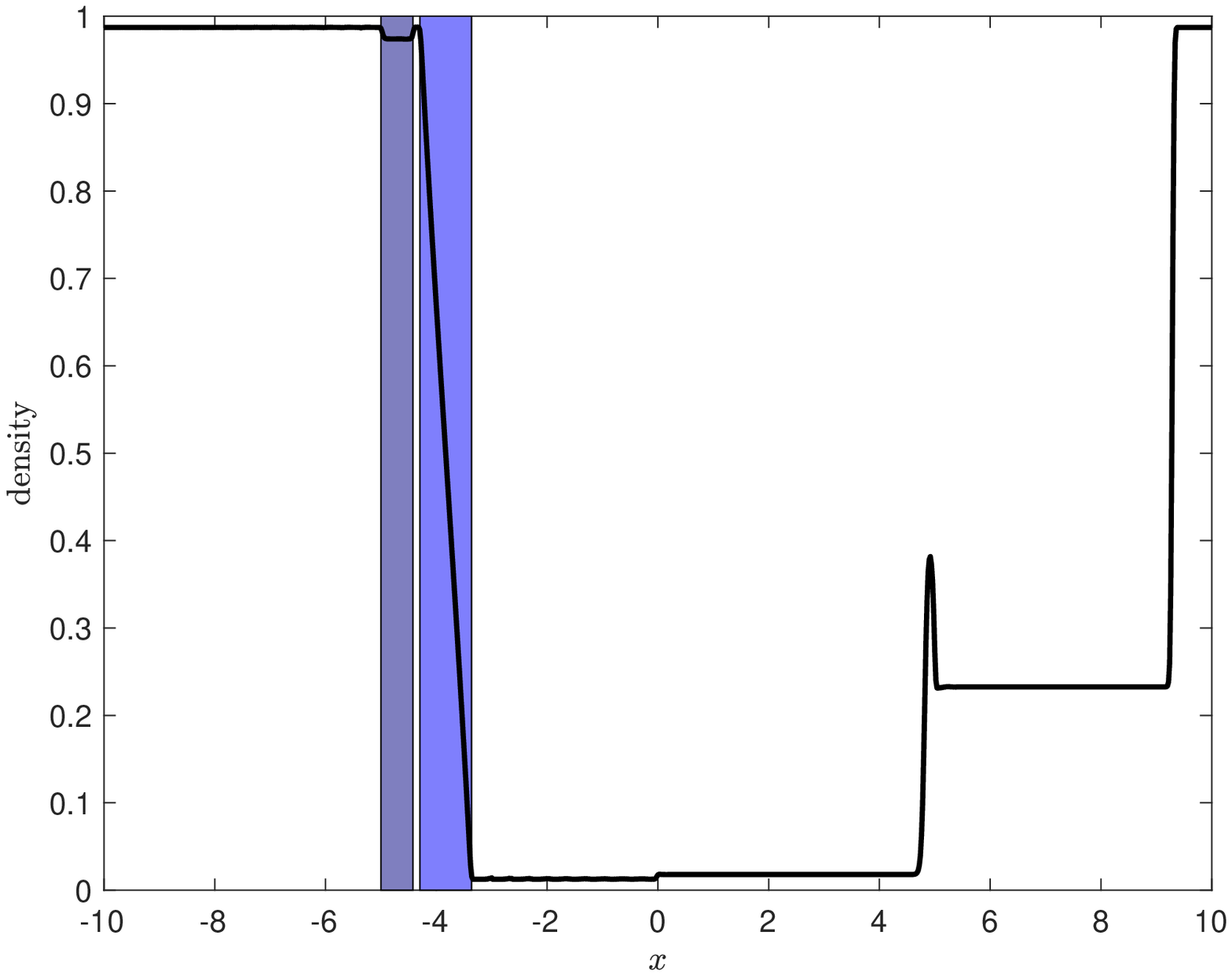}
\label{fig:Ex2amplePathBeta0d}
}
\end{figure}

\begin{figure}[H]
\subfigure[$t = 26.8$: second and third accident removed.]{
\includegraphics[width = 0.45\textwidth]{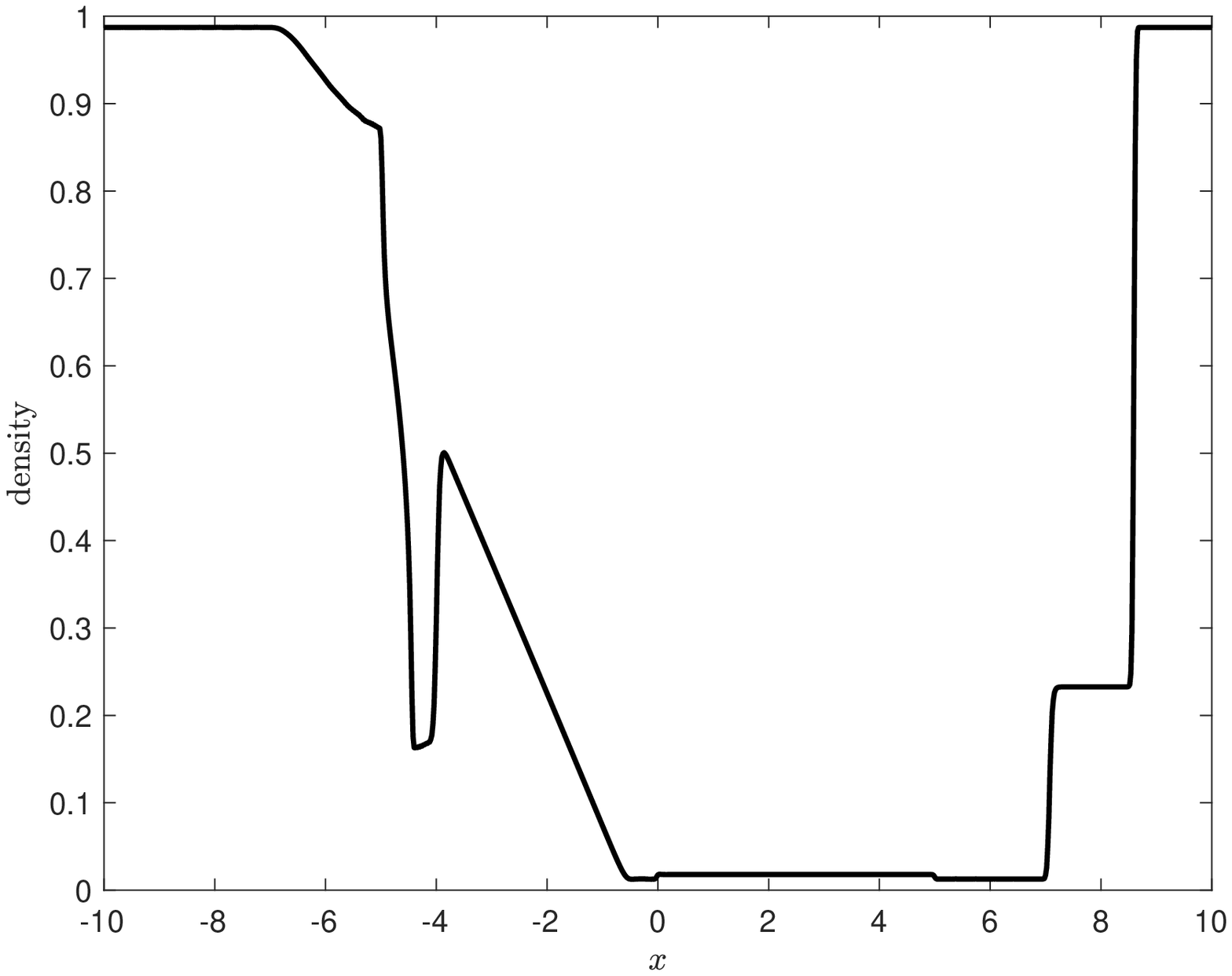}
\label{fig:Ex2SamplePathBeta0e}
}
\hfill
\subfigure[$t = 29.35$: fourth accident.]{
\includegraphics[width = 0.45\textwidth]{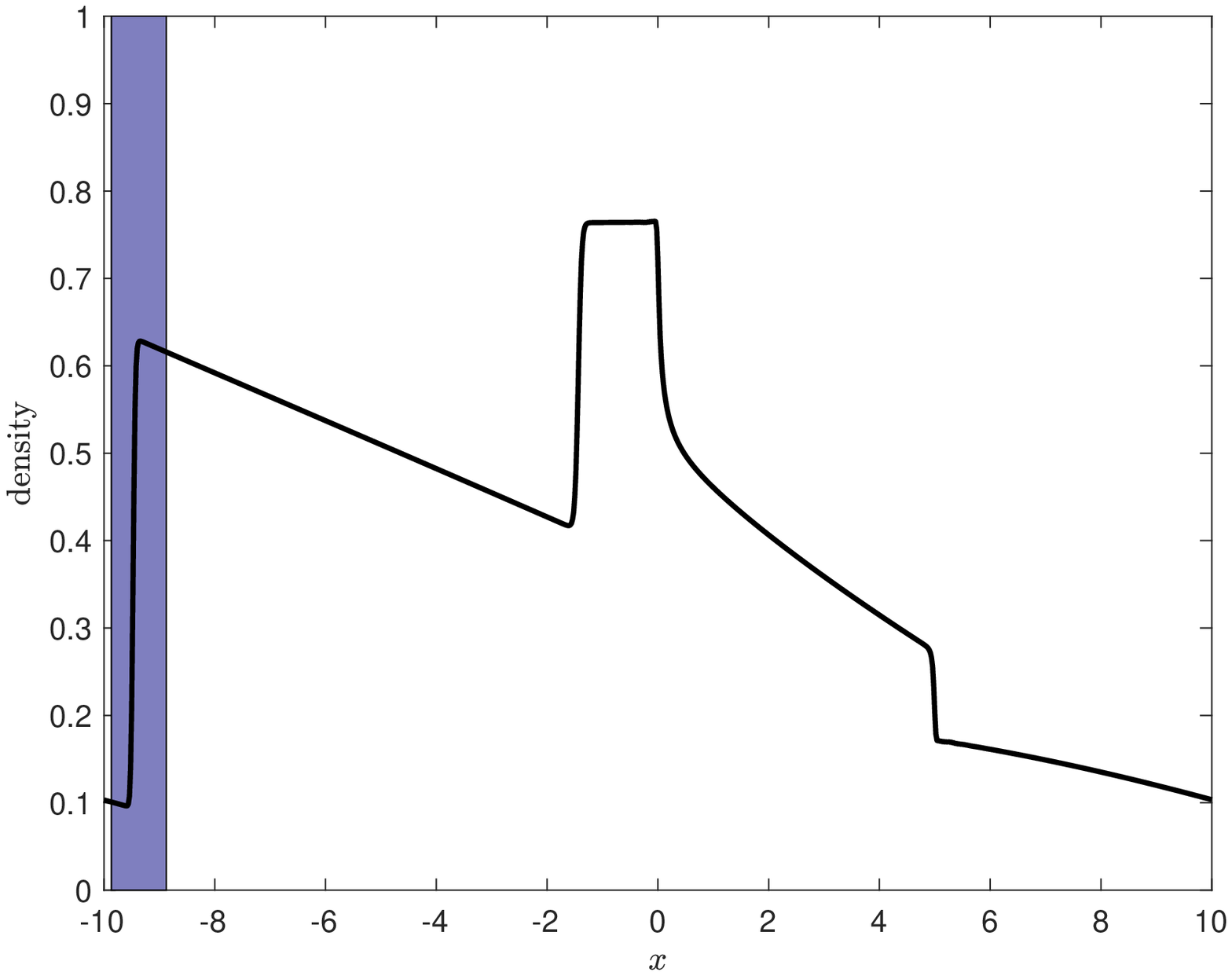}
\label{fig:Ex2SamplePathBeta0f}
}
\caption{$\beta = 0$.}
\label{fig:Ex2SamplePathBeta0}
\end{figure}

Figure \ref{fig:Ex2SamplePathBeta0} shows the traffic density (black bold line) for different points in time and using only the information of $D\rho^+$ to determine the position of an accident, i.e.\ $\beta = 0$. The rectangles in the figures indicate the range of the road affected by an accident, where a bright color corresponds to a capacity drop of $0.99$ and the other color of $0.5$, see $\mu^c$ in  \eqref{eq:muEx1}. Since the initial distribution is constant with a value of $0.4$, we draw the density at the first time at which an accident happens in Figure \ref{fig:Ex2SamplePathBeta0a}. Due to a spatial inhomogeneous road capacity $c_{road}(x)$, the initial density profile changed to a non-constant equilibrium traffic density. As we would expect, the accident happens at the incresaing part of the density, i.e.\ at the end of the traffic jam, with a road capacity reduction of $0.99$. At this position a traffic jam occurs until the accident is removed, see Figure \ref{fig:Ex2SamplePathBeta0b}.  The traffic density relaxes to an equilibrium density again and the second accident happens at the end of the traffic jam as Figure \ref{fig:Ex2SamplePathBeta0c} indicates. Again a capacity reduction of $0.99$ has been randomly chosen and a third accident occurs right after the second accident. The latter can be seen in Figure \ref{fig:Ex2amplePathBeta0d}, which shows the traffic density at the time, where the second accident gets resolved. 

\begin{figure}[H]
\subfigure[$t = 7.4$: first accident.]{
\includegraphics[width = 0.45\textwidth]{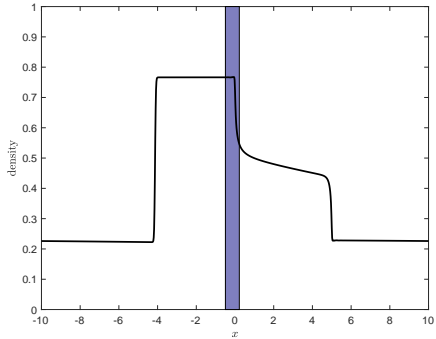}
\label{fig:Ex1SamplePathBeta05a}
}
\hfill
\subfigure[$t = 8.35$: second accident.]{
\includegraphics[width = 0.45\textwidth]{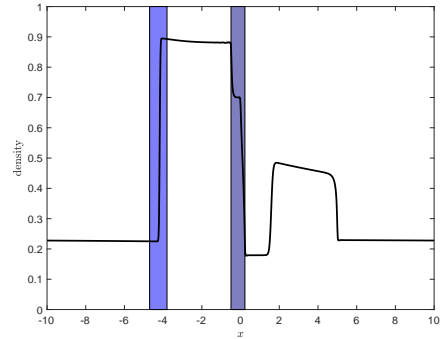}
\label{fig:Ex1SamplePathBeta05b}
}
\subfigure[$t = 9.95$: third accident and first accident removed.]{
\includegraphics[width = 0.45\textwidth]{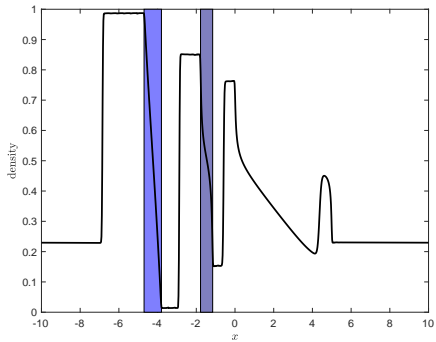}
\label{fig:Ex1SamplePathBeta05c}
}
\hfill
\subfigure[$t = 10.5$: fourth accident within second accident.]{
\includegraphics[width = 0.45\textwidth]{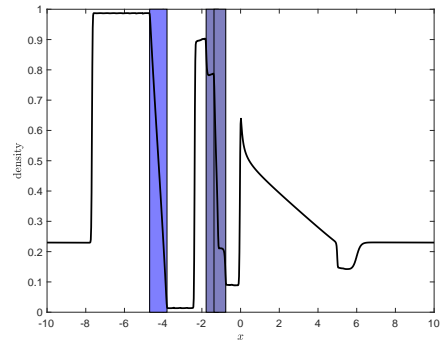}
\label{fig:Ex1SamplePathBeta05d}
}
\caption{$\beta = 0.5$.}
\label{fig:Ex1SamplePathBeta05}
\end{figure}

At the time, where both accidents are resolved, see Figure \ref{fig:Ex2SamplePathBeta0e}, we see the high impact of the previous accidents on the density, which does not reach the equilibrium state until the next accident occurs as Figure \ref{fig:Ex2SamplePathBeta0f} shows. Again, the position of the accident is at an increasing part of the density.
Altogether, we see that our model is able to map the ideas of accidents at places with an increasing density and the numerical solutions look very confident using the CFL condition with equality. 
In the following, we discuss simulation results using the parameter $\beta = 0.5$ shown in Figure \ref{fig:Ex1SamplePathBeta05}. We face an approximately equilibrium density at the time of a first accident again, see Figure \ref{fig:Ex1SamplePathBeta05a}. Here, the accident occurs close to the position zero, which is not an increasing part of the density. The accident is therefore created by the flux, which is uniform on the interval [-10,10] while the density is close to equilibrium. 

As Figure \ref{fig:Ex1SamplePathBeta05b} shows, the second accident happens at the traffic jam end. After the first accident has been removed, a third accident occurs and Figure \ref{fig:Ex1SamplePathBeta05c} shows the traffic density at the time right before the fourth accident occurs. The fourth accident is inside the area of the second accident and has a small size of impact, see Figure \ref{fig:Ex1SamplePathBeta05d}. The latter accident occurred at this position since the flux around $\rho = 0.5$ is the most highest and we are not in a stationary state.

\begin{figure}[H]
\subfigure[$\beta = 0$.]{
\includegraphics[width = 0.45\textwidth]{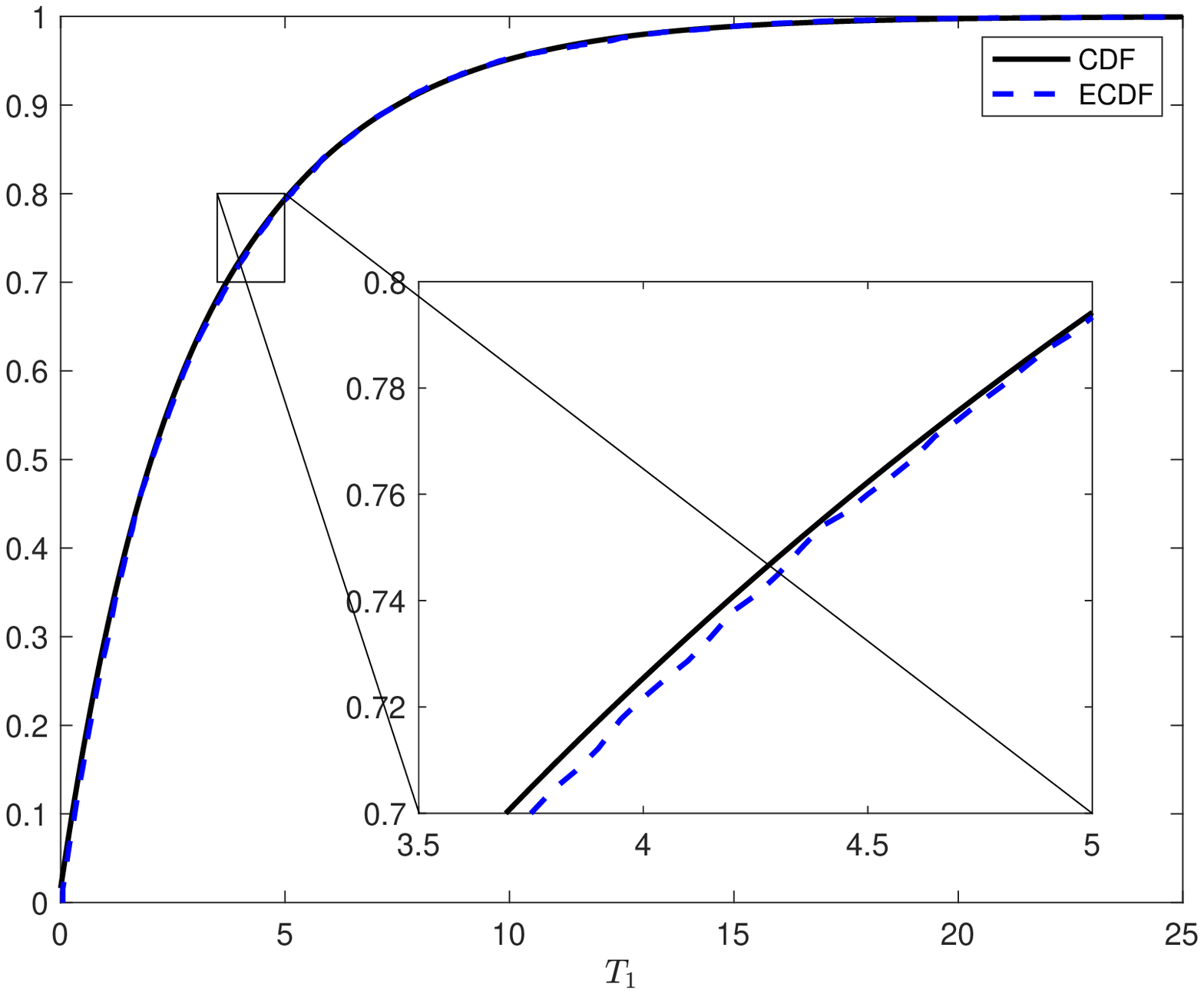}
\label{fig:ECDFvsCDFBeta0}
}
\hfill
\subfigure[$\beta = 0.5$.]{
\includegraphics[width = 0.45\textwidth]{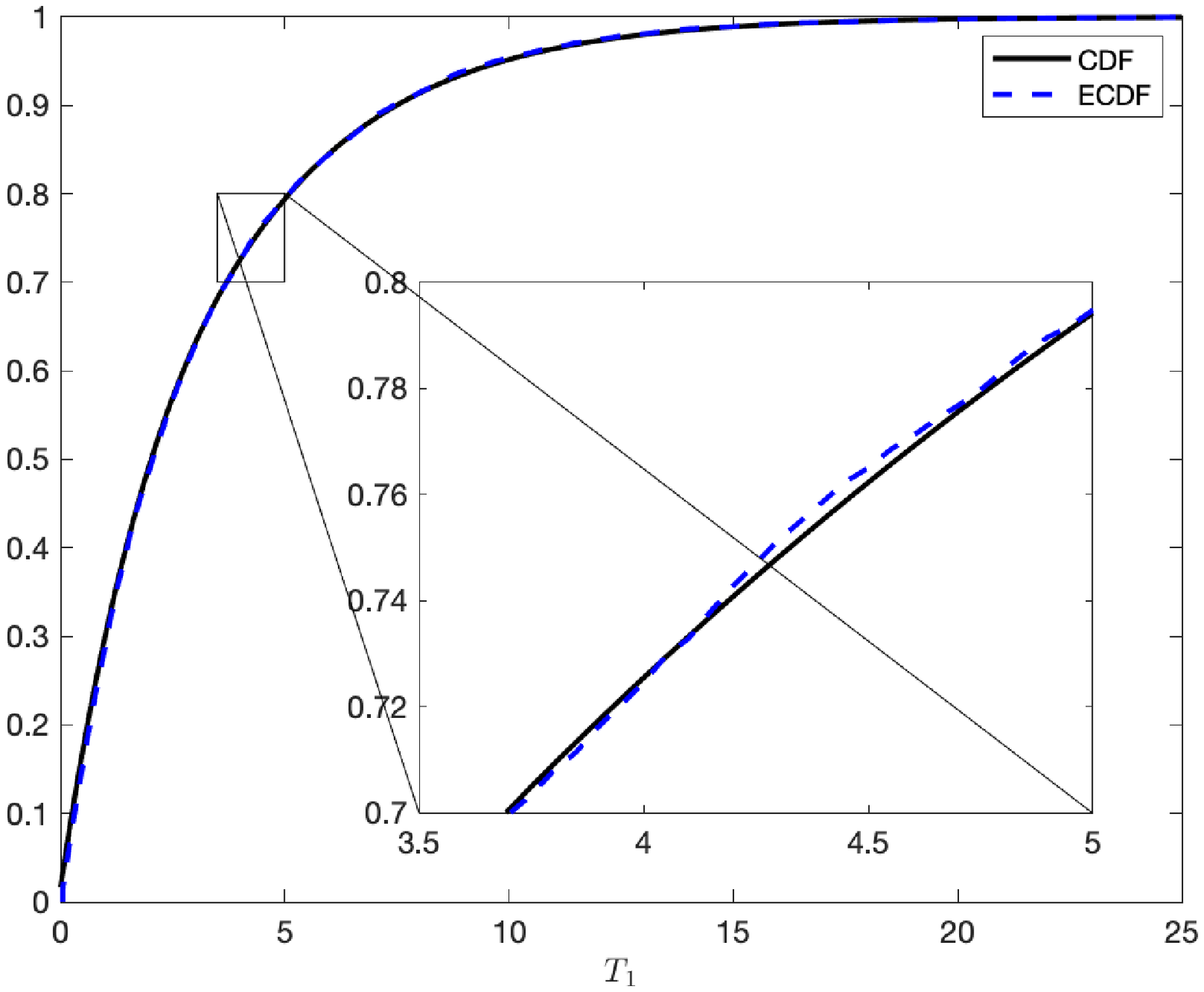}
\label{fig:ECDFvsCDFBeta05}
}
\caption{ECDF of the first accident time $T_1^a$ compared with the CDF  for $T_1$ in \eqref{eq:ThinningJumpDist}.}
\label{fig:ECDFvsCDF}
\end{figure}

\textit{Numerical verification of the approximate scheme.}
In order to verify numerically that the approximate algorithm works well, we study the distribution of the first jump time, i.e.\ the first time of an accident. Formula \eqref{eq:ThinningJumpDist} exactly describes the cumulative distribution function (CDF), which can be approximated using the Lax-Friedrichs scheme to approximate $\phi$.  Using the left-sided rectangular rule to approximate $\int_0^t \psi(\phi_\tau(y_0))d\tau$ and the Matlab function \texttt{ecdf}\footnote{Documentation: https://de.mathworks.com/help/stats/ecdf.html, 2019} to compute the empirical cumulative distribution function (ECDF) yields the results shown in Figure \ref{fig:ECDFvsCDF} computed by using $10^4$ samples of the first accident time $T_1^a$.
First of all, we observe a  very good fitting of the CDF by the ECDF computed with the approximation Algorithm \ref{alg:approx_jumps}. This implies that the corresponding probability distributions are close (in the weak sense). Furthermore, we observe that the parameter $\beta$ has no significant influence on the shape or values of the CDF as Figures \ref{fig:ECDFvsCDFBeta0} and \ref{fig:ECDFvsCDFBeta05} show. 

\begin{figure}[H]
\subfigure[$\beta = 0$.]{
\includegraphics[width = 0.45\textwidth]{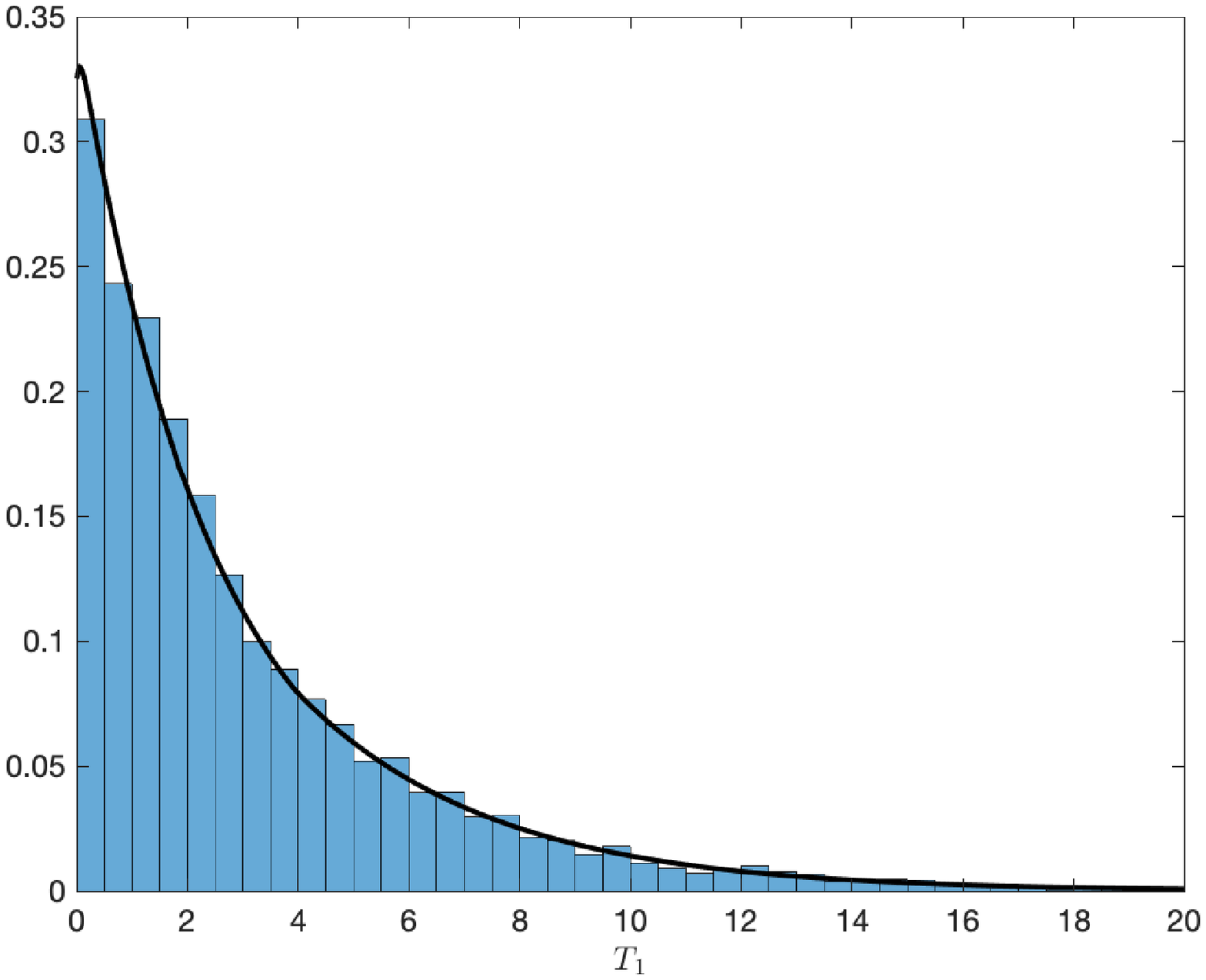}
\label{fig:HistogramBeta0}
}
\hfill
\subfigure[$\beta = 0.5$.]{
\includegraphics[width = 0.45\textwidth]{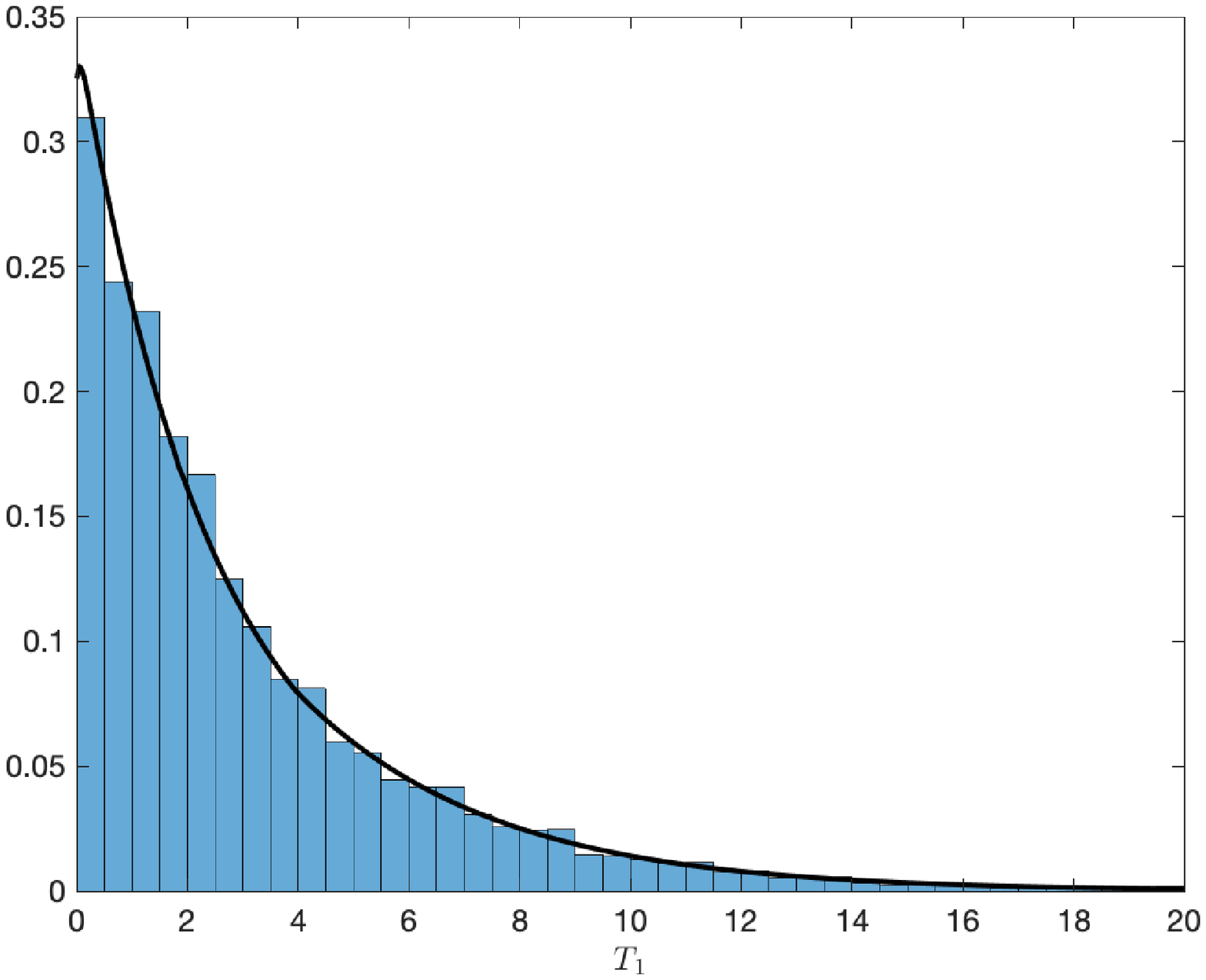}
\label{fig:HistogramBeta05}
}
\centering
\subfigure{
\includegraphics[width = 0.45\textwidth]{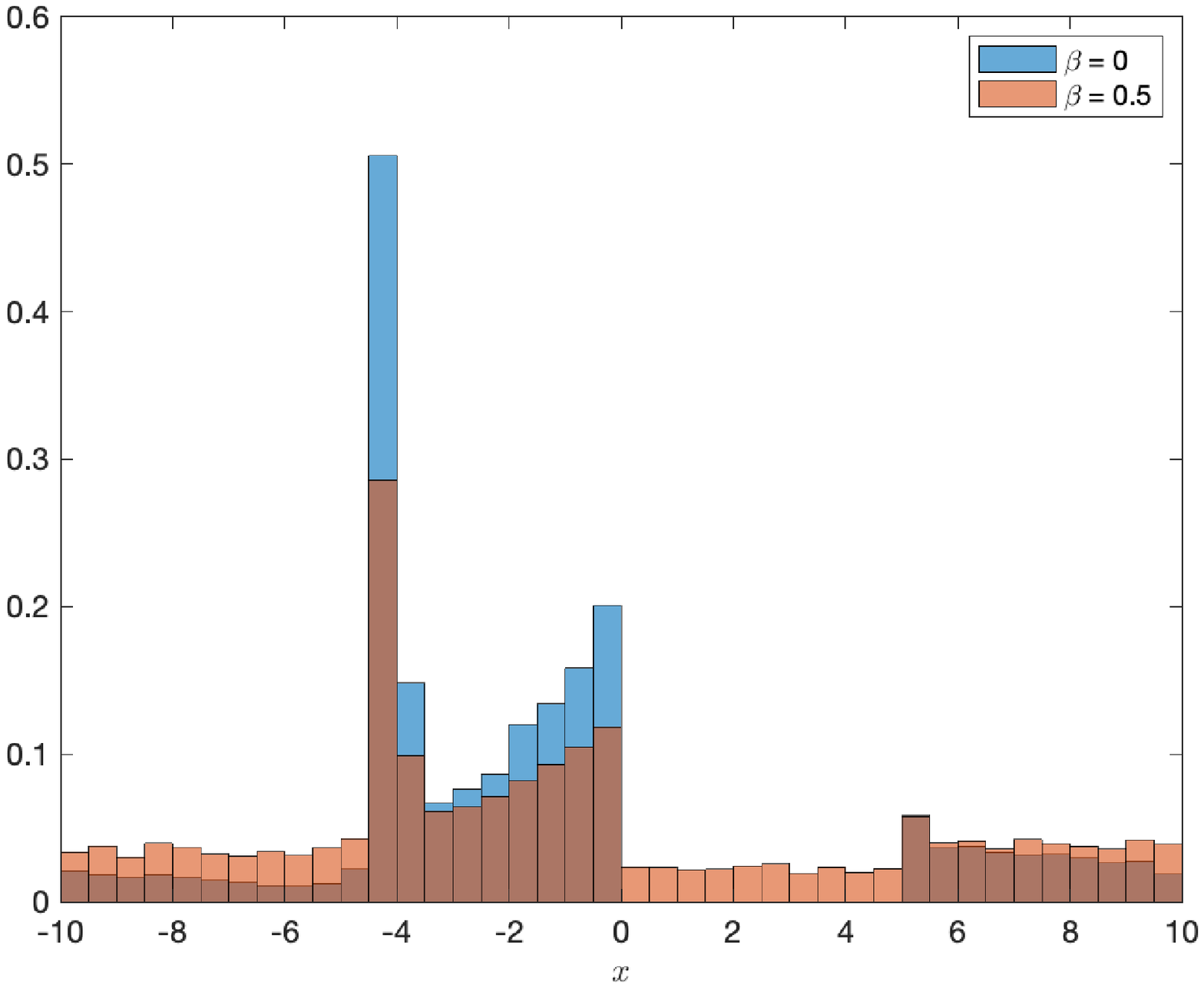}
}
\caption{Histograms of $T^a_1$ (first row)
and of the first accident's position (second row).}
\label{fig:HistogramPositions}
\end{figure}

In order to compare a histogram generated by the approximation procedure with the exact probability density function (pdf) $g(t)$, we can differentiate  \eqref{eq:ThinningJumpDist} and obtain
\begin{align*}
g(t) = \psi(\phi_t(y_0))e^{-\int_0^t \psi(\phi_\tau(y_0))d\tau}.
\end{align*}
Figure \ref{fig:HistogramPositions} shows a histogram of samples of $T_1^a$ and the theoretical result $g(t)$. We observe a good agreement between both quantities again, also independent of the choice of $\beta$. 

Finally, we discuss the distribution of the first accident's position. Figure \ref{fig:HistogramPositions} shows the histogram of samples of the first accident's position, where we distinguish the cases $\beta = 0$ and $\beta = 0.5$ again. In both cases, the probability having an accident at position $x = -4$ is the most highest, which corresponds to the congestion end in the stationary traffic profile, see Figure \ref{fig:Ex2SamplePathBeta0c} for example. One significant difference between $\beta = 0$ and $\beta = 0.5$ can be observed for $x \in [0,5]$, where in the case of $\beta = 0$, i.e.\ no flux information, no accident happens. 

In contrast, for $\beta =0.5$, there is a strictly positive probability having an accident in this interval, which is clear since the stationary value of $\rho$ is approximately at the maximal flow, i.e.\ at $0.5$.

To conclude, the numerical simulations inherit the ideas for the stochastic traffic flow model and the numerical results are convincing.

\section{Conclusion}
We successfully have derived a stochastic traffic flow model capturing random traffic accidents. Furthermore, a tailored numerical approximation scheme has been introduced, which also has been validated in numerical simulation examples. 

The stochastic traffic flow model allows for road capacity planning and controlling variable speed limit systems in such a way that traffic accidents are rarely events, which might be future research. Additionally, the extension to a second order traffic models and networks can be considered.

\section*{Acknowledgments}
This work was supported by the BMBF project ENets (05M18VMA) and the DFG project GO1920/10-1.
\bibliographystyle{siam}
\bibliography{references}
\end{document}